\documentclass[12pt]{article}

\usepackage{amsmath,amssymb,amsfonts,amsthm,graphics,verbatim}

\include{diagram}

 \DeclareMathOperator{\nd}{Nbhd}

\DeclareMathOperator{\rk}{\text{rank}}

\DeclareMathOperator{\sln}{\bold{SL_n}}

\DeclareMathOperator{\ar}{\bold{G}(\mathcal{O}_S)}
\DeclareMathOperator{\HOM}{\widetilde{H}}
\DeclareMathOperator{\gd}{\bold{G_{\derived}}}
\DeclareMathOperator{\derived}{der}

\DeclareMathOperator{\AP}{\bold{A}}

\DeclareMathOperator{\torus}{\bold{T}}
\DeclareMathOperator{\cent}{\bold{Z}}

\begin{document}

\newcommand{\cdottie}{}
\newcommand{\si}{\mathcal{O}_S}
\newcommand{\ka}{\kappa}
\newcommand{\ga}{\Gamma}
\newcommand{\tga}{\Gamma}
\newcommand{\tg}{G}
\newcommand{\rt}{\rightarrow}
\newcommand{\gmga}{\tga \backslash \tg}
\newcommand{\neu}{N(\Gamma)}
\newcommand{\ap}{{\mathcal{A}}}
\newcommand{\lb}{\lambda}
\newcommand{\se}{\subseteq}
\newcommand{\e}{\varepsilon}
\newcommand{\nre}{\sqrt[n]{\varepsilon}}
\newcommand{\tp}{\widetilde{\phi}}
\newcommand{\pp}{\partial \phi}
\newcommand{\app}{{\mathcal{A}}^+}
\newcommand{\ppr}{\partial \phi _R}
\newcommand{\pps}{{\partial \phi _R}_*}
\newcommand{\pu}{U _{\partial}}
\newcommand{\gh}{\widehat {\Gamma}}
\newcommand{\ghp}{\widehat {\Gamma} ^{\varphi}}
\newcommand{\hp}{H^+}
\newcommand{\oppr}{\overline{\partial \phi _R}}
\newcommand{\ta}{K \gamma_\mathfrak{p}}
\newcommand{\de}{\delta}
\newcommand{\tgd}{K \gamma_\mathfrak{p}(\de)}
\newcommand{\p}{\phi}
\newcommand{\pn}{\phi ^{-1}}
\newcommand{\g}{\gamma}
\newcommand{\pgd}{\pi(\gr ,\de)}
\newcommand{\ep}{\hfill $\blacksquare$ \\ \medskip}
\newcommand{\pru}{\bigskip \noindent {\bf Proof:} \;}
\newcommand{\sa}{\measuredangle _}
\newcommand{\xy}{X_\infty}
\newcommand{\gy}{\gamma_\infty}
\newcommand{\py}{\pi_\infty}
\newcommand{\ey}{e_\infty}
\newcommand{\xr}{X_{\mathfrak{p}}}
\newcommand{\gr}{\gamma_{\mathfrak{p}}}
\newcommand{\pr}{\pi_{\mathfrak{p}}}
\newcommand{\er}{e_{\mathfrak{p}}}
\newcommand{\phr}{\phi_{\mathfrak{p}}}
\newcommand{\phy}{\phi_{\infty}}

\input{picinpar.sty}

\title{Finiteness properties of arithmetic groups over function fields}

\author{Kai-Uwe Bux and Kevin Wortman\thanks{Supported by an NSF Postdoctoral
Fellowship.}}

\maketitle

\begin{abstract}
We determine when an arithmetic subgroup of a reductive group
defined over a global function field is of type $FP_\infty$ by
comparing its large-scale geometry to the large-scale geometry of
lattices in real semisimple Lie groups.
\end{abstract}

\bigskip \noindent \Large \textbf{1 Introduction} \normalsize \bigskip

Throughout this paper, $K$ is a global function field, and $S$ is
a finite nonempty set of pairwise inequivalent valuations on $K$.
We let $\mathcal{O}_S \leq K$ be the corresponding ring of
$S$-integers. We denote a reductive $K$-group by $\mathbf{G}$.

In 1971 Serre proved that $\mathbf{G}(\si)$ is of type $WFL$ if
and only if the semisimple $K$-rank of $\mathbf{G}$ equals $0$;
see Th\'{e}or\`{e}me 4 of \cite{Serre initial} and the following
Compl\'{e}ments.

As type $FP_\infty$ is a weaker property than type $WFL$, an
immediate consequence is that $\mathbf{G}(\si)$ is of type
$FP_\infty$ if the semisimple $K$-rank of $\mathbf{G}$ equals $0$.
 The converse of this statement had
been believed since the late 1970's\footnote{see e.g. the final
introductory paragraph of \cite{Stuhler}} and evidence had been
collected to support it as a conjecture. However, it remained
unresolved in general.

Our main result confirms this conjecture:

\bigskip \noindent \textbf{Theorem A} \emph{The arithmetic group
$\ar$ is of type $FP_{\infty}$ if and only if the semisimple
$K$-rank of $\mathbf{G}$ equals $0$.}

\bigskip

As a special case of our main result,
$\mathbf{SL_n}(\mathbb{F}_q[t])$ is not of type $FP_\infty$. Even
this basic example was previously unknown in full generality; see
Example below. Here, $\mathbb{F}_q[t]$ is the ring of polynomials
with one indeterminate $t$ and coefficients in the finite field
with $q$ elements, $\mathbb{F}_q$.

We will also give a more precise statement about the finiteness
lengths of arithmetic groups; see Theorem B. As a special case of
that result, $\mathbf{SL_n}(\mathbb{F}_q[t])$ is not even of type
$FP_{n-1}$.

\smallskip \noindent \textbf{Historical remarks.} Interest in
 the finiteness properties of arithmetic groups
  over function fields was sparked
in 1959 by Nagao's proof that $\mathbf{SL_2}(\mathbb{F}_q[t])$ is
not finitely generated \cite{Nagao}.

Activities of the next 33 years completely determined which
arithmetic subgroups of reductive groups
  over function fields are finitely generated,
and which are finitely presented (the answers fit the form of
Conjecture C below). Work on these results was carried out by
Behr, Hurrelbrink, Keller, Kneser, Lubotzky, McHardy, Nagao,
O'Meara, Rehmann-Soul\'{e}, Serre, Splitthoff, and Stuhler. See
\cite{Behr reduction}, \cite{Behr example}, \cite{Behr survey},
\cite{Hurrelbrink}, \cite{Keller}, \cite{Lubotzky finite
generation}, \cite{Mc}, \cite{O'Meara}, \cite{R-S}, \cite{Serre},
\cite{Splitthoff}, and \cite{Stuhler 1}.

Less understood are the higher finiteness properties for these
groups, such as type $FP_n$ for $n \geq 3$. Aside from the result
of Serre mentioned earlier, all of the work in this direction has
been carried out with heavy restrictions on $\mathbf{G}$ and
$\mathcal{O}_S$; see the papers of Abels, Abramenko, Behr, and
Stuhler  (\cite{Abels}, \cite{Abramenko thesis}, \cite{Abramenko
paper}, \cite{Abramenko book}, \cite{Behr preprint}, and
\cite{Stuhler}).

\bigskip

Theorem A follows as a corollary of Theorem B below. Before
presenting the statement of Theorem B, we introduce some notation.

\smallskip \noindent \textbf{Type $\mathbf{FP_m}$.} Recall that for a
commutative ring $R$, we say a group $\Gamma$ is \emph{of type
}$FP_m$ \emph{over} $R$ if there exists a projective resolution
$$P_m \rt P_{m-1} \rt \cdot \cdot \cdot \rt P_1 \rt P_0 \rt R \rt 0$$
of finitely generated $R \Gamma$ modules, where the action of $R
\Gamma$ on $R$ is trivial. If $\Gamma$ is of type $FP_m$ over $R$
for all nonnegative integers $m$, we say that $\Gamma$ is \emph{of
type} $FP_\infty$ \emph{over} $R$. If $\Gamma$ is of type $FP_m$
(resp. $FP_\infty$) over $\mathbb{Z}$, we simply write that
$\Gamma$ is \emph{of type }$FP_m$ (resp. $FP_\infty$).

\smallskip \noindent \textbf{Remark.} Every group is of type
$FP_0$. Type $FP_1$ is equivalent to the property of finite
generation. Every finitely presented group is of type $FP_2$, but
the converse does not hold in general; see Bestvina-Brady's
Example 6.3(3) in \cite{Best-Brad}.

\smallskip \noindent \textbf{Finiteness length.} The
\emph{homological finiteness length of} $\Gamma$ \emph{over} $R$
is defined to be the nonnegative integer
$$\phi (\Gamma ;R)=\sup \{ m \mid \Gamma \text{ is of type } FP_m \text{ over } R\}$$
For short, we write $\phi (\Gamma )$ in place of $\phi (\Gamma
;\mathbb{Z})$.

\smallskip \noindent \textbf{Rank.} For any field extension $L/K$, the $L$-\emph{rank}
of $\mathbf{G}$, denoted $\rk _L \mathbf{G}$, is the dimension of
a maximal $L$-split torus of $\mathbf{G}$. The \emph{semisimple}
$L$-\emph{rank} of $\mathbf{G}$ is the $L$-rank of the derived
subgroup of $\mathbf{G}$.

If $\rk _L (\mathbf{G}) =0$, we say $\mathbf{G}$ is
$L$-anisotropic. Otherwise, $\mathbf{G}$ is $L$-isotropic.

\smallskip \noindent \textbf{Sum of local ranks.} For a valuation
$v$ of $K$, let $K_v$ be the completion of $K$ with respect to
$v$. For any $K$-group $\mathbf{G}$, we define the nonnegative
integer
$$k(\mathbf{G},S)=  \sum _{v \in S} \rk _{K_v}
\mathbf{G} $$ We are now prepared to state

\bigskip \noindent \textbf{Theorem B} \emph{If $\mathbf{H}$ is a connected noncommutative
absolutely almost
simple $K$-isotropic $K$-group, then}
$$\phi \big(\mathbf{H}(\si ) \big) \leq k(\mathbf{H},S)-1$$
That Theorem A follows from Theorem B is routine; see e.g. 2.6(c)
of \cite{Behr survey}.

\bigskip

\noindent \textbf{Example.} A special case of Theorem B is the
inequality $$\phi \Big(\sln (\mathbb{F}_q[t])\Big) \leq n-2$$ or
more generally,
$$\phi \Big(\sln (\mathcal{O}_S)\Big) \leq |S|(n-1)-1$$ Indeed,
for any field $L$, the number $\rk _L \sln$ equals the dimension
of the diagonal subgroup in $\sln$. Hence, for any $K$ and any
$S$, we have
$$k(\sln , S)=\sum _{v \in S} \rk _{K_v} \sln=|S|(n-1)$$

This inequality is known to be sharp in some cases. For example,
Stuhler showed that $\phi (\mathbf{SL_2} (\mathcal{O}_S))=|S|-1$
\cite{Stuhler}, and Abels and Abramenko independently showed that
$\phi (\mathbf{SL_n} (\mathbb{F}_q[t]))=n-2$ as long as $q \geq
2^{n-2}$ or $q \geq \tbinom{n-2}{\lfloor \frac{n-2}{2} \rfloor }$
respectively \cite{Abels}, \cite{Abramenko thesis}.

\smallskip \noindent \textbf{Is the inequality sharp in general?} Theorem B provides
evidence for the following long-standing conjecture, which offers
a striking relation between the two functions $\phi$ and $k$.

\bigskip \noindent \textbf{Conjecture C} \emph{If $\mathbf{H}$ is a connected noncommutative
absolutely almost simple $K$-isotropic $K$-group, then}
$$\phi \big(\mathbf{H}(\si ) \big) = k(\mathbf{H},S)-1$$
See \cite{Behr survey} for other evidence.

\bigskip

\smallskip \noindent \textbf{Type $\mathbf{F_m}$.} Recall that a group
$\Gamma$ is of type $F_m$ if there exists an Eilenberg-Mac\,Lane
complex $K(\Gamma , 1)$ with finite $m$-skeleton. For $m \geq 2$,
a group is of type $F_m$ if and only if it is finitely presented
and of type $FP_m$. It then follows from \cite{Behr survey} that
$FP_m$ and $F_m$ are equivalent conditions
 for the arithmetic groups considered in this paper. Thus,
  Theorems A and B, and Conjecture C, may be equivalently stated
  by substituting $F_m$ for $FP_m$.

\smallskip \noindent \textbf{Type $\mathbf{WFL}$.} Although we
will make no further use of it, we recall the definition of type
$WFL$ for completeness with respect to comments in the initial
portion of the introduction: A group $\Gamma$ is \emph{of type}
$WFL$ if there exists a torsion-free finite-index subgroup of
$\Gamma$, and if for any such subgroup $\Gamma '$, the ring
$\mathbb{Z}$ admits a finite length resolution by finitely
generated free $\mathbb{Z}\Gamma '$-modules.

\smallskip \noindent \textbf{Contrast with number fields.} Our theorems are
particular to the case of global fields of positive
characteristic. In characteristic zero, we have the following

\bigskip \noindent \textbf{Theorem D (Raghunathan, Borel-Serre)} \emph{Any $S$-arithmetic
 subgroup of a reductive
group defined over a global number field contains a finite-index
torsion-free subgroup $\Gamma$ that allows for a finite $K(\Gamma
, 1)$.}

\bigskip \noindent In particular, any $S$-arithmetic subgroup of a reductive
group in characteristic zero is finitely presented and of type
$FP_\infty$. Examples of groups for which the above theorem
applies include $\mathbf{SL_n}(\mathbb{Z})$ and
$\mathbf{SL_n}(\mathbb{Z}[1/p])$. See \cite{Raghunathan paper} for
the case of arithmetic groups, and \cite{Borel-Serre} for the case
of $S$-arithmetic groups.

\smallskip \noindent \textbf{Idea behind the proof.} Although
Theorem A shows a difference between arithmetic groups in positive
characteristic and those in characteristic zero, it is through the
deep-rooted similarity of these two families that we shall find a
proof of Theorem A.

Indeed, our motivating example for proving Theorem B was the proof
of Epstein-Thurston that $\mathbf{SL_3}(\mathbb{Z})$ is not
combable; see Chapter 10 Section 4 of \cite{word processing}.
Recall that their proof proceeds by creating an exponential Dehn
function for $\mathbf{SL_3}(\mathbb{Z})$ as follows. A family of
closed curves with increasing lengths are constructed in a portion
of the symmetric space
$\mathbf{SL_3}(\mathbb{R})/\mathbf{SO_3}(\mathbb{R})$ that is a
bounded distance from the subset $\mathbf{SL_3}(\mathbb{Z})
\mathbf{SO_3}(\mathbb{R}) \se
\mathbf{SL_3}(\mathbb{R})/\mathbf{SO_3}(\mathbb{R})$. The discs in
$\mathbf{SL_3}(\mathbb{R})/\mathbf{SO_3}(\mathbb{R})$ that fill
these loops in the most metrically efficient manner have areas
that are quadratic in the length of the loops that they bound.
These efficient discs are not so useful in studying the
large-scale geometry of $\mathbf{SL_3}(\mathbb{Z})$ though, since
the discs travel farther away from the subspace
$\mathbf{SL_3}(\mathbb{Z})
 \mathbf{SO_3}(\mathbb{R})$ as the length of their boundary curves
increase.

To better understand the geometry of $\mathbf{SL_3}(\mathbb{Z})$,
we only consider filling discs of the constructed loops that are
contained in the original bounded neighborhood of
$\mathbf{SL_3}(\mathbb{Z}) \mathbf{SO_3}(\mathbb{R})$. What is
shown in \cite{word processing} is that any such family of discs
would be metrically inefficient in the sense that the discs would
have areas that are at least exponential in the length of their
boundary curves. The result is an exponential Dehn function for
$\mathbf{SL_3}(\mathbb{Z})$, which implies that
$\mathbf{SL_3}(\mathbb{Z})$ is not combable.

Our proof of Theorem B in the special case when $\mathbf{H}(\si)
=\mathbf{SL_3}(\mathbb{F}_q[t])$ proceeds by constructing an
analogous family of loops in a bounded neighborhood of a given
$\mathbf{SL_3}(\mathbb{F}_q[t])$-orbit in the Euclidean building,
$X$, associated to $\mathbf{SL_3}\big (\mathbb{F}_q((t^{-1}))
\big)$ where $\mathbb{F}_q((t^{-1}))$ is a field of formal Laurent
series. As with the case for $\mathbf{SL_3}(\mathbb{Z})$, the
closed curves have metrically efficient filling discs in $X$ whose
areas are quadratic in the length of their boundary curves. Also
like the case for $\mathbf{SL_3}(\mathbb{Z})$, these filling discs
travel farther away from the given orbit as the length of their
boundary curves increase, so they are not helpful in learning
about the large-scale geometry of
$\mathbf{SL_3}(\mathbb{F}_q[t])$. However, in this case, there
does not exist a filling disc for any of our constructed loops
that is contained in the fixed bounded neighborhood of the
$\mathbf{SL_3}(\mathbb{F}_q[t])$-orbit. Indeed, $X$ is a
contractible $2$-dimensional simplicial complex, so filling discs
are essentially unique. We then apply K. Brown's filtration
criterion to conclude that $\mathbf{SL_3}(\mathbb{F}_q[t])$ is not
finitely presented.

\smallskip \noindent \textbf{Distortion dimension.} The
contrast between arithmetic groups over function fields with
arithmetic groups over number fields diminishes if we consider a
metric analogue of finiteness length.

Let us direct our attention for the moment to an irreducible
lattice $\Gamma$ in a semisimple group over arbitrary nondiscrete
locally compact fields; we can even allow for $\Gamma$ to be
nonarithmetic.  We let $X_\Gamma$ be the natural product of
irreducible symmetric spaces and Euclidean buildings that $\Gamma$
acts on. Given a point $x\in X_\Gamma$ and a real number $r$, we
define the space
$$X_\Gamma(r)
 = \{\, y \in  X_\Gamma \mid d(y\,,\,\Gamma x) \leq r \,\}$$

Using the Hurewicz theorem---as in Abels-Tiemeyer's Theorem 1.1.4
of \cite{A-T}---and recalling that type $FP_m$ and type $F_m$ are
equivalent conditions for $\Gamma$ allows us to state K. Brown's
filtration criterion for $\Gamma$ to be of type $FP_m$ in terms of
homotopy groups. Precisely, $\Gamma$ is of type $FP_m$ if and only
if for any real number $r \geq 0$ there exists a real number $r'
\geq r$ such that for any $k<m$ the homomorphism induced by
inclusion
$$\pi _k  \big(X_\Gamma (r) \,,\, x \big) \longrightarrow
 \pi _k \big(X_\Gamma (r') \,,\, x \big)$$ is trivial; see
 Theorems 2.2 and 3.2 \cite{Brown filtration}.

 Hence, if $\Gamma$ acts cocompactly on
$X_\Gamma$, it is of type $FP_\infty$. If $\Gamma$ does not act
cocompactly, then Theorems A and D (along with Corollary 7.3 of
\cite{Lubotzky finite generation}) characterize those $\Gamma$
contained in semisimple groups over function fields as precisely
those which fail to be of type $FP_\infty$. (Recall that an
arithmetic lattice $\mathbf{H}(\si)$ acts cocompactly on
$X_{\mathbf{H}(\si)}$ if and only if the absolutely almost simple
$K$-group $\mathbf{H}$ is $K$-anisotropic.)

To include metric properties of the large-scale geometry of
lattices, we define $\Gamma$ as being \emph{undistorted up to
dimension} $m$ if: given any $r \geq 0$, there exist real numbers
$r' \geq r$, $K \geq 1$, and $C \geq 0$ such that for any $k<m$
and any Lipschitz $k$-sphere $s \se X_\Gamma (r)$, there exists a
Lipschitz $(k+1)$-ball
    $B _\Gamma \se X_\Gamma(r')$ with $\partial B _\Gamma = s$ and
    $$\text{volume}(B _\Gamma ) \leq K [\text{volume} (B_X)] +C$$ for
    all Lipschitz $(k+1)$-balls
    $B_X \se X$ with $\partial B_X = s$.
We adopt the convention that $\Gamma$ is always undistorted up to
dimension $0$.

Now we define the \emph{distortion dimension of} $\Gamma$ to be
the nonnegative integer $$\psi (\Gamma )=\sup \{ m \mid \Gamma
\text{ is undistorted up to dimension } m \}$$

\bigskip \noindent \textbf{Conjecture E} \emph{If $\Gamma$ is an irreducible lattice
in a semisimple group over nondiscrete locally compact fields,
then $\psi \big( \Gamma \big)=\infty$ if and only if $\Gamma $
acts cocompactly on $X_\Gamma$.}

\bigskip

That $\psi (\Gamma )=\infty$ when $\Gamma $ acts cocompactly is
clear. The converse had been conjectured for lattices in real
semisimple Lie groups following the Epstein-Thurston proof that
$\psi (\mathbf{SL_n}(\mathbb{Z}))\leq n-2$, and a general proof
seems approachable. (See 10.4 \cite{word processing} for
$\mathbf{SL_n}(\mathbb{Z})$.)

Less attention has been given to $S$-arithmetic lattices in
characteristic zero, but the conjecture should not change in this
setting. In positive characteristic, Conjecture E follows from
Theorem A.

As Conjecture E extends
 Theorem A for absolutely almost simple $K$-groups $\mathbf{G}$ into
 the context of arbitrary global fields, we are naturally
    led to speculate how Conjecture C might be broadened
    to include fields of characteristic zero.
     Thus, we define $\tau (X_\Gamma)$ to be the Euclidean rank of $X_\Gamma$, and
we note that for $\Gamma = \mathbf{H}(\si )$ as in Conjecture C,
$\tau (X_\Gamma)= k (\mathbf{H},S)$. We ask

\bigskip \noindent \textbf{Question F} \emph{Let $\Gamma$ be a noncocompact irreducible lattice
in a semisimple group over nondiscrete locally compact fields. Is
it true that }
$$\psi \big( \Gamma \big)=\tau (X_\Gamma ) -1$$ \emph{If not, then can the definition of $\psi$
be reasonably modified so that the above formula is true? }

\bigskip

This problem is daunting. For example, an affirmative answer to
the first question implies Thurston's claim that
$\mathbf{SL_n}(\mathbb{Z})$ has a quadratic Dehn function for $n
\geq 4$.
 See also remarks from Gromov's book (5.D.(5).(c)
\cite{Gromov}).

We will just mention a few pieces of evidence for a positive
answer. We note that
    $\psi (\Gamma )=0$ if and only if either
    $\Gamma$ is not finitely generated or the word metric on
    $\Gamma$ is not quasi-isometric to the metric induced from
    its action on $X_\Gamma$; see the example in
    Section 2. Hence, it follows from work of Lubotzky
(Corollary 7.3 of \cite{Lubotzky finite generation}) and
Lubotzky-Mozes-Raghunathan \cite{L-M-R} that $\psi (\Gamma )=0$ if
and only if $\tau (X_\Gamma ) =1$. Note also that Leuzinger-Pittet
\cite{Leuzinger-Pittet}, Behr \cite{Behr survey}, and (a
generalization of) Taback (Lemma 4.2 of \cite{Taback}), show that
$\tau (X_\Gamma ) =2$ implies $\psi (\Gamma )=1$.

For related material, see the papers of Dru\c{t}u, Hattori,
Leuzinger-Pittet, Noskov, and Pittet: \cite{Drutu 1}, \cite{Drutu
2}, \cite{Hattori}, \cite{Leuzinger-Pittet 2}, \cite{Noskov}, and
\cite{Pittet}.

\smallskip \noindent \textbf{Possible generalizations to other rings of functions.} In
\cite{B-W} we use techniques from this paper to give a geometric
proof that $\mathbf{SL_2}(\mathbb{Z}[t,t^{-1}])$
 is not finitely presented---a
fact first proved by Krsti\'{c}-McCool \cite{K-M}.

It is likely that the ideas below can be used to do more in this
direction of generalizing Theorem B to apply to a class of fields
and rings that properly includes global function fields and their
rings of $S$-integers; see the question in the introduction of
\cite{B-W}.

\smallskip \noindent \textbf{Outline of the paper.} We begin in
Section 2 with a special case of our proof to motivate what
follows. The proof of Theorem B is contained in Section 3. In the
appendix we include the proofs of two well-known results for
completeness: the existence of anisotropic maximal tori in
semisimple groups defined over nondiscrete locally compact fields
of positive characteristic, and the ``if" implication of Theorem
A.

\bigskip

\noindent \textbf{Acknowledgements.} We thank Stephen DeBacker and
Dan Barbasch for telling us about the existence of anisotropic
tori in semisimple groups over local fields of positive
characteristic.

We are also happy to thank Ross Geoghegan for suggesting that our
proof would be more efficiently carried out by using homology with
coefficients in $\mathbb{F}_p$, and Indira Chatterji for helping
us to see the connection between finite generation for lattices in
nonarchimedean semisimple groups and bounded distortion of the
word metric for lattices in real semisimple Lie groups.

This paper also benefitted from conversations we had with our
colleagues at Cornell University and the University of Chicago:
 Tara Brendle, Nathan Broaddus, Kenneth Brown, Allen Hatcher, Paul
Jung, Robert Kottwitz, Alessandra Pantano, Ravi Ramakrishna,
Edward Swartz, Karen Vogtmann, and Dan Zaffran; it is our pleasure
to thank them.

The second author thanks Benson Farb and Dan Margalit for their
instrumental encouragement.

\bigskip \noindent \Large \textbf{2 An example} \normalsize \bigskip

The first piece of evidence for Theorem A was:

\bigskip \noindent \textbf{Theorem 2(a) (Nagao)} \emph{The group
$\mathbf{SL_2}(\mathbb{F}_q[t])$ is not finitely generated.}

\bigskip

In this section we will see how our proof of Theorem B applies to
this special case. For motivation, we will first review some of
the geometry of $\mathbf{SL_2}(\mathbb{Z})$, a mathematical cousin
of $\mathbf{SL_2}(\mathbb{F}_q[t])$.

 \begin{window}[0,r,{\InputIfFileExists{figure_1.tex}{}{}},{}]
      Consider the action of
      $\mathbf{SL_2}(\mathbb{Z})$ on the hyperbolic plane $\mathbb{H}^2$.
      The diagram shows the upper half-plane
      model. There is a distinguished point
      $\infty$ at the top of the diagram
       that no $\mathbf{SL_2}(\mathbb{Z})$-orbit accumulates on.
      Specifically, it is well-known that the orbit of the complex number
      $i$ avoids the open horoball
      $B$ that is centered at
      $\infty$ and consists of all
complex numbers with imaginary parts greater than $1$. The
boundary of this horoball is approximated by the points $n+i$ for
$n\in \mathbb{Z}$. (Notice that $n+i={{1 \; n} \choose {0 \; 1} }
 i$.)
  \end{window}

The geodesic joining $i$ and $1+i$ travels into the horoball $B$.
The geodesic between $i$ and $2+i$ travels farther into the
horoball, the geodesic between $i$ and $3+i$ farther still, and so
on. Continuing this process, we see that no metric neighborhood of
the orbit $\mathbf{SL_2}(\mathbb{Z})  i\se \mathbb{H}^2$ is convex
in $\mathbb{H}^2$. Sufficiently large metric neighborhoods of
$\mathbf{SL_2}(\mathbb{Z})  i$ are however connected, as
$\mathbf{SL_2}(\mathbb{Z})$ is finitely generated:

    \bigskip \noindent \textbf{Lemma 2(a)} \emph{Suppose a finitely generated group $\Gamma$ acts on a geodesic
metric
      space $X$. Then, for any point $x \in X$, there is
      a number $R>0$ such that the $R$-neighborhood}
       $$ \nd_R(\Gamma   x) \se X$$
      \emph{is connected.}

    \pru
      Let $\{\gamma_1,\gamma_2,\ldots,\gamma_n\}$ be a finite generating set
      for $\Gamma$. Choose $R$ such that the
      ball $B_R(x)$ contains all translates $\gamma_i    x$.
      Then $\Gamma    B_R(x)
      =\nd_R(\Gamma    x)$
      is connected.

    \ep

Now let's look at a proof of Nagao's theorem along the lines of
our proof for Theorem B. This is not Nagao's original proof,
rather it is a simplified form of Stuhler's argument
\cite{Stuhler}.

The locally compact field $\mathbb{F}_q((t^{-1}))$ contains
$\mathbb{F}_q[t]$ as a discrete subring. Thus,
$\mathbf{SL_2}(\mathbb{F}_q[t])$ is a discrete subgroup of
$\mathbf{SL_2}\big(\mathbb{F}_q((t^{-1}))\big)$.

There is a natural nonpositively curved space that
$\mathbf{SL_2}\big(\mathbb{F}_q((t^{-1}))\big)$ acts on: the
      regular $(q+1)$-valent tree
      $T$. If $\mathbb{F}_q[[t^{-1}]] < \mathbb{F}_q((t^{-1}))$ is
      the ring of Taylor series, then this well-known action is
      obtained by identifying the vertices of $T$ with homothety classes
      of spanning
$\mathbb{F}_q[[t^{-1}]]$-submodules of a $2$-dimensional vector
space over $\mathbb{F}_q((t^{-1}))$ that are free and of rank $2$.
This is in analogy to the identification of the unit tangent
bundle of $\mathbb{H}^2$ with the unit tangent bundle of the
Teichm\"{u}ller space of $2$-dimensional Euclidean tori with
volume equal to $1$, or equivalently, with homothety classes of
spanning $\mathbb{Z}$-submodules of $\mathbb{R}^2$ that are free
and of rank $2$.

      Just as the boundary of $\mathbb{H}^2$ is a circle, or $\mathbb{P}^1(\mathbb{R})$,
      the boundary of $ T$ can be identified with
$\mathbb{P}^1\big(\mathbb{F}_q((t^{-1})) \big)$. We use the
standard identification of
$\mathbb{P}^1\big(\mathbb{F}_q((t^{-1})) \big)$ with
$\mathbb{F}_q((t^{-1})) \cup \{\infty \}$. The group
$\mathbf{SL_2}\big(\mathbb{F}_q((t^{-1}))\big)$ has two induced
actions on $\mathbb{P}^1\big(\mathbb{F}_q((t^{-1})) \big)$: one
from its action on $T$, and one from its action on the
$2$-dimensional vector space $\mathbb{F}_q((t^{-1})) ^2$. These
actions coincide.

\begin{window}[0,r,{\InputIfFileExists{figure_2.tex}{}{}},{}]
In another analogy with the situation for
$\mathbf{SL_2}(\mathbb{Z})$, any ray from any point $x \in T$
towards $\infty$ escapes every metric neighborhood of the orbit
$\mathbf{SL_2}(\mathbb{F}_q[t])  x$ (one can see this using
Mahler's compactness criterion). The diagram on the right has $x$
contained in the geodesic joining the two boundary points $0$ and
$\infty $. We write $f * x \in T$ as shorthand for the point ${{1
\; f} \choose {0 \; 1} }  x$. The geodesic segment between $x$ and
$t^n * x$ is the portion of the geodesic joining the boundary
points $0$ and $t^n$ that lies at or above the level of $x$ in the
diagram. These segments contain increasing subsets of the geodesic
ray from $x$ to $\infty$ as $n \in \mathbb{N}$ grows. Hence, no
metric neighborhood of the orbit $\mathbf{SL_2}(\mathbb{F}_q[t])
 x\se T$ is convex.
\end{window}

The comparison with $\mathbf{SL_2}(\mathbb{Z})$ stops here since,
in $T$, convexity is equivalent to connectedness. Using Lemma
2(a), we see that $\mathbf{SL_2}( \mathbb{F}_q[t])$ is not
finitely generated. Our proof is complete.

\bigskip \noindent \Large \textbf{3 Proof of Theorem B} \normalsize \bigskip

In what follows, we let $p$ equal the characteristic of $K$.
Rather than proving Theorem B directly, our goal will be to prove
a slightly stronger claim:

  \bigskip \noindent \textbf{Proposition B*} \emph{If $\mathbf{H}$ is a
  connected noncommutative absolutely almost
simple $K$-isotropic $K$-group, then}
  $$\phi \big( \mathbf{H}(\mathcal{O}_S) ; \mathbb{F}_p \big) \leq k(\mathbf{H},S)-1$$

\bigskip

    Theorem B follows since a group $\Gamma$ is of type $FP_m$ over $\mathbb{F}_p$ if it is
    of type $FP_m$ over $\mathbb{Z}$: just tensor a projective resolution for
    $\mathbb{Z}$ by finitely generated
 $\mathbb{Z}\Gamma$-modules with $\mathbb{F}_p$ to obtain
    a projective resolution for $\mathbb{F}_p$ by finitely generated
    $\mathbb{F}_p \Gamma$-modules.

\medskip \noindent \large \textbf{3.1 Method of proof for Proposition
B*} \normalsize \medskip

We define the ring $$K_S=\prod _{v\in S} K_v $$ so that
$$\mathbf{H}(K_S)=\prod _{v\in S} \mathbf{H}(K_v)$$ Let $X$ be the
Euclidean building corresponding to $\mathbf{H}(K_S)$, that is the
product of the irreducible Euclidean buildings for
$\mathbf{H}(K_v)$. Recall that $X$ has dimension $k(\mathbf{H},
S)$.

We fix a base point $e \in X$ (to be specified later) and
    consider closed metric neighborhoods of the orbit $\mathbf{H}(\si
    )  e$. That is, for each number $r \geq 0$, we set
    $$ X( r)
      =
      \big\{
        x \in X
        \mid
        d\big( x\,,\,\mathbf{H}(\si )  e \big)
        \leq
        r
      \big\}$$
    We will find a number $r_0>0$ and construct, for each
    $r \geq r_0$, a cycle in $ X(r_0)$
    that represents a nontrivial element in the reduced homology
    group with coefficients in $\mathbb{F}_p$ $$\HOM_{k(\mathbf{H},S)-1}\big(X(
    r)\,;\,\mathbb{F}_p\big)$$
This shows that the inclusions
    $X(r_0)\subseteq X(r)$
    induce nontrivial homomorphisms $$ \HOM_{k(\mathbf{H},S)-1}\big(X(r_0)\,;\,\mathbb{F}_p\big)
      \longrightarrow
      \HOM_{k(\mathbf{H},S)-1}\big( X(r)\,;\,\mathbb{F}_p \big) $$

In view of K. Brown's filtration criterion (see Theorem 2.2 and
the following remark in \cite{Brown filtration}), the existence of
this family of nontrivial homomorphisms together with the
following standard facts about the action of $\mathbf{H}(\si )$ on
$X$ implies Proposition B*: (\emph{i}) $X$ is contractible;
(\emph{ii}) $\mathbf{H}(\si)$ acts on $X$ with finite cell
stabilizers; and (\emph{iii}) the subspaces $X(r)$ are
$\mathbf{H}(\si)$-invariant and compact modulo $\mathbf{H}(\si)$.

\smallskip \noindent \textbf{Excluding a tree.} For the remainder of this paper, we
will assume that $k(\mathbf{H},S) > 1$. That is, we assume that
$X$ is not a tree. This assumption is made only to avoid
complications in our exposition; the philosophy of the proof still
applies to the case when $k(\mathbf{H},S)=1$ as is shown in
Section 2.

\medskip \noindent \large \textbf{3.2 An apartment
 coarsely separated by $\mathbf{H}\mathbf{(\si)}$} \normalsize \medskip

We will find an apartment in $X$ that ``coarsely intersects" an
$\mathbf{H}(\mathcal{O}_S)$-orbit in a hyperplane. (Later, we will
use this $(k(\mathbf{H},S)-1)$-dimensional hyperplane and its
translates to construct the $(k(\mathbf{H},S)-1)$-cycles mentioned
above.) Since apartments in $X$ correspond to products of maximal
$K_v$-split tori in $\mathbf{H}$, this problem reduces to algebra.

We begin by choosing a parabolic group that will follow with us
throughout our proof. In what follows, we are assuming that the
reader has a basic knowledge of the structure of parabolic
subgroups of reductive groups relative to fields that are not
algebraically closed, as can be found for example in 21.11 and
21.12 of \cite{Borel text}.

Since $\mathbf{H}$ is $K$-isotropic, there exists a nontrivial
maximal $K$-split torus of $\mathbf{H}$. We let $\Phi _K$ be the
roots of $\mathbf{H}$ with respect to this torus. Choose an
ordering on $\Phi _K$, and let $\Delta _K \se \Phi_K$ denote the
corresponding collection of simple roots.

Choose, and fix throughout, a root $\alpha _ 0 \in \Delta _K$. We
define the $1$-dimensional $K$-split torus $$\mathbf{T_1}=\Big(
\bigcap _{\alpha \in \Delta _K - \alpha _ 0
}\text{ker}(\alpha)\Big) ^\circ$$ The above superscript $\circ$
denotes the connected component of the identity. We let
$\mathbf{Z_H(T_1)}$ be the centralizer of $\mathbf{T_1}$ in
$\mathbf{H}$.

There exists a maximal proper $K$-parabolic subgroup of
$\mathbf{H}$, denoted $\mathbf{P^+}$, with the following Levi
decomposition:
$$\mathbf{P^+}=\mathbf{R_u(P^+)} \rtimes \mathbf{Z_H(T_1)}$$
In the above, $\mathbf{R_u(P^+)}$ is the unipotent radical of
$\mathbf{P^+}$.

We can expand the Levi decomposition to a Langlands decomposition
by noting that $\mathbf{Z_H(T_1)}$ is an almost direct product of
$\mathbf{T_1}$, the derived group
$\mathbf{Z_H(T_1)}_{\text{der}}$, and $\mathbf{D_a}$ for some
$K$-anisotropic diagonalizable group $\mathbf{D_a}$. Thus:
$$\mathbf{P^+}=\mathbf{R_u(P^+)} \rtimes \mathbf{T_1D_a}\mathbf{Z_H(T_1)}_{\text{der}}$$

Before proceeding with the existence of the torus and the
apartment that is our goal in this section, we record the
following well-known result.

\bigskip \noindent \textbf{Proposition 3.2(a)} \emph{Let
$\mathbf{G}$ be a reductive $K$-group. Then for any finite
nonempty set $S'$ of pairwise inequivalent valuations and any
family $\{ \mathbf{A}_v \}_{v\in S'}$ of maximal $K_v$-tori of
$\mathbf{G}$, there is a maximal $K$-torus $\mathbf{A}_\pi$ of
$\mathbf{G}$ and group elements $g_v \in \mathbf{G}(K_v)$ such
that
$$\mathbf{A}_\pi=\sideset{^{g_v}}{_v}{\AP}$$
for all $v \in S'$, where $\sideset{^{g_v}}{}{\AP}_v$ denotes $\AP
_v$ conjugated by $g_v$. }

\pru There is a proof of this proposition in Section 7.1 Corollary
3 of \cite{Pl-Ra} for the case when $K$ is a global number field.
The proof also applies for global function fields after replacing
the argument for the $K$-rationality of the variety of maximal
tori in $\mathbf{G}$ with the proof of Theorem 7.9 in
\cite{Bo-Sp}.

\ep

We will make use of the above proposition in the proof of the
proposition below.

    \bigskip \noindent \textbf{Proposition 3.2(b)} \emph{There exists a maximal
     $K$-torus $\mathbf{A} \leq \mathbf{H}$
such that:}
      \begin{itemize}
        \item[(\emph{i})]
          \emph{The maximal $K$-split torus of $\mathbf{A}$ is $\mathbf{T_1}$, and}
        \item[(\emph{ii})]
          $\mathbf{A}$ \emph{contains a maximal $K_v$-split torus of} $\mathbf{H}$ \emph{for all} $v\in S$.
      \end{itemize}

    \pru For each $v \in S$,
let $\mathbf{A}_v$
    be a maximal $K_v$-torus of $\mathbf{Z_H(T_1)}_{\text{der}}$ such that $\mathbf{A}_v$ contains
     a maximal $K_v$-split
    torus of $\mathbf{Z_H(T_1)}_{\text{der}}$.

    Then choose a valuation
    of $K$, call it $w$, that is inequivalent to any of the valuations of $S$,
     and let $\mathbf{A}_w$ be a maximal
    $K_w$-torus in $\mathbf{Z_H(T_1)}_{\text{der}}$ that is
    $K_w$-anisotropic. The existence of such a torus is well-known;
     see the appendix for a proof.

Now apply Proposition 3.2(a) to
$\mathbf{G}=\mathbf{Z_H(T_1)}_{\text{der}}$ and $S'=S\cup \{w\}$.
Since $\mathbf{A}_\pi$ is $K_w$-anisotropic, it is necessarily
$K$-anisotropic.
 Therefore, part (\emph{i}) is satisfied by
$$\mathbf{A}= \mathbf{T_1}  \mathbf{D_a}^\circ \mathbf{A}_\pi$$

To verify part (\emph{ii}), note that $\mathbf{T_1}$ is contained
in a maximal $K_v$-split torus of $\mathbf{H}$. Hence,
\begin{align*} \rk _{K_v}
(\mathbf{H} ) & = \rk _{K_v} (\mathbf{Z_H(T_1)} ) \\
& = \rk _{K_v} (\mathbf{T_1D_a} ) + \rk _{K_v} (\mathbf{Z_H(T_1)}_{\text{der}}) \\
& = \rk _{K_v} (\mathbf{T_1D_a} ) + \rk _{K_v}
 (\mathbf{A}_\pi ) \\
 & = \rk _{K_v} (\mathbf{A} )
\end{align*}
\enlargethispage*{33pt}
 \ep

Since $\mathbf{A}$ contains a maximal $K_v$-split torus for all
$v\in S$, there is an apartment $\Sigma  \se X$ that
$\mathbf{A}(K_S)$ acts on properly and cocompactly as a
translation group of maximal rank, $\text{dim}(\Sigma )=
k(\mathbf{H},S)$. By Dirichlet's units theorem (see Theorem 5.12
\cite{Pl-Ra}) and the preceding proposition, the arithmetic group
$\mathbf{A}(\si )$ is a finitely generated abelian group of rank
$$ \Big(\sum _{v \in S} \rk _{K_v} (\mathbf{A})\Big) - \rk _K
(\mathbf{A}) = k(\mathbf{H},S) -1 $$

Choose a point $e \in \Sigma $. Since $\mathbf{A}(\si ) \leq
\mathbf{A}(K_S )$ acts properly on $\Sigma$, the base point $e$ is
contained in an affine hyperplane $V \se \Sigma$, of dimension
$k(\mathbf{H},S)-1$, that $\mathbf{A}(\si)$ acts on cocompactly.
This point $e \in \Sigma$ is the point we specify for our
definition in Section 3.1 of the spaces $X(r)\se X$.

\smallskip \noindent \textbf{Example.} In the case when $K=\mathbb{F}_q(t)$,
$\si =\mathbb{F}_q[t]$, and $\mathbf{H}=\mathbf{SL_3}$, the torus
$\mathbf{T_1}$ can be taken as the group of matrices of the form
\[
\begin{pmatrix}
a  & 0 & 0 \\
0 & a & 0 \\
0 & 0 & a^{-2}
\end{pmatrix}
\]
Then the parabolic group $\mathbf{P^+}$ can be taken to be the
determinate $1$ matrices of the form \[
\begin{pmatrix}
* & * & * \\
* & * & * \\
0 & 0 & *
\end{pmatrix}
\]
The resulting group $\mathbf{R_u(P^+)}$ would be the
$2$-dimensional commutative group
\[\begin{pmatrix}
1 & 0 & * \\
0 & 1 & * \\
0 & 0 & 1
\end{pmatrix}
\]
This would leave the semisimple group
$\mathbf{Z_H(T_1)}_{\text{der}}$ to be the copy of $\mathbf{SL_2}$
that sits in the upper left corner of $\mathbf{SL_3}$, and
$\mathbf{D_a}$ would be trivial.

The group $\mathbf{A}(\si )$ in this example can be taken to be
the group generated by the matrix \[
\begin{pmatrix}
t^2 +1  & t & 0 \\
t & 1 & 0 \\
0 & 0 & 1
\end{pmatrix}
\]
With the notation from the proof of Proposition 3.2(b), the
Zariski closure of $\mathbf{A}(\si)$ would equal $\mathbf{A}_\pi
\leq \mathbf{Z_H(T_1)}_{\text{der}}$.

\medskip \noindent \large \textbf{3.3 A space for
 the manufacture of cycles: choosing $\mathbf{r_0}$} \normalsize \medskip

 Let $D\se V$ be a fundamental
domain for the $\mathbf{A}(\si )$-action on $V$. By Behr-Harder
reduction theory (e.g. Satz 3 of \cite{Behr reduction}), there is
a compact set $C \se \mathbf{R_u(P^+)}(K_S)$ such that $
\mathbf{R_u(P^+)}(K_S) = \mathbf{R_u(P^+)} (\si ) C$. Since
$\mathbf{A}(\si )$ normalizes $\mathbf{R_u(P^+)}(K_S)$, we have:

    \begin{eqnarray*}
      \mathbf{R_u(P^+)}(K_S)  V
      &\subseteq&
      \mathbf{R_u(P^+)}(K_S) \mathbf{A}(\mathcal{O}_S)  D \\
      &\subseteq&
      \mathbf{A}(\mathcal{O}_S)\mathbf{R_u(P^+)}(K_S)  D \\
      &\subseteq&
      \mathbf{A}(\mathcal{O}_S)\mathbf{R_u(P^+)}(\mathcal{O}_S)C  D \\
      &\subseteq&
      \mathbf{H}(\mathcal{O}_S)C  D
    \end{eqnarray*}

    Since the region $C  D \se X$ is bounded, we can choose a number $r_0>0$ such that
     $$ \mathbf{R_u(P^+)}(K_S)  V
      \subseteq
      \mathbf{H}(\mathcal{O}_S)C  D
      \subseteq
      X(r_0) $$

It is inside the space $ \mathbf{R_u(P^+)}(K_S)  V$ where we shall
produce cycles that remain nontrivial in the homology of $X(r)$
for $r \geq r_0$.

\medskip \noindent \large \textbf{3.4 A direction away from
$\mathbf{X(r_0)}$} \normalsize \medskip

Recall our choice of $\alpha _ 0 \in \Delta _K$ from the beginning
of Section 3.2. This root is nontrivial when restricted to
$\mathbf{T_1}$, so for any $v \in S$, the set $$\{\, a \in
\mathbf{T_1} (K_v) \mid |\alpha _ 0(a)|_v>1 \,\}$$ is nonempty and
open in the Hausdorff topology induced by the metric $|\cdot |_v$
on $K_v$ that arises from $v$.

Since $\mathbf{T_1}$ is $K$-isomorphic to an affine line with a
point removed, it satisfies the weak approximation property with
respect to $S$. That is, the diagonal embedding
$$\mathbf{T_1}(K) \longrightarrow \prod _{v\in S} \mathbf{T_1}(K_v)$$
has a dense image. Therefore, there exists some $a_+ \in
\mathbf{T_1}(K)$ such that $$|\alpha _ 0 (a_+)|_v >1$$ for all
$v\in S$.

It will be important for us later to have a direction in $\Sigma$
that leads away from every $X(r)$. The direction we will use is
given by the sequence $(a_+^n  e)_{n \in \mathbb{N}}$. Note that
the above condition on $a_+$ assures us that the sequence
$(a_+^n)_{n \in \mathbb{N}}$ is not contained in any compact
subset of $\mathbf{H}(K_S)$. Therefore, $(a_+ ^n  e)_{n \in
\mathbb{N}}$ does specify a direction in $\Sigma$.

Let $ X ^\infty$ be the visual boundary of $X$. It can be
identified in a natural way with the spherical Tits building for
$\mathbf{H}(K_S)$. Note that $X^\infty$ is the spherical join of
the spherical buildings for the groups $\mathbf{H}(K_v)$ with $v
\in S$. We let $\Sigma ^\infty \se X^\infty$ be the apartment
corresponding to $\Sigma \se X$, and we let $a_+^\infty \in \Sigma
^\infty$ be the accumulation point of $(a_+^n  e)_{n \in
\mathbb{N}}$.

We let $\Pi _+^\infty$ be
      the unique simplex in $X^\infty$ that is maximal among all simplices stabilized by the
      action of $\mathbf{P^+}(K_S)$. Note
      that $\Pi _+^\infty$ is the spherical join over $S$ of the simplices
      associated with
       $\mathbf{P^+}(K_v)$ in
      the spherical buildings for $\mathbf{H}(K_v)$.

The conditions on our choice of $a_+$ were imposed to insure that
the orbit of $e$ under its iterates would accumulate inside $\Pi
_+ ^\infty$. Specifically, we have:

\bigskip \noindent \textbf{Lemma 3.4(a)} \emph{The point $a_+^\infty \in \Sigma ^\infty$
is contained in $\Pi _+ ^\infty$.}

\pru Using the definition of spherical joins, we can reduce to the
case when $S$ contains a single valuation $v$. What follows is
routine; see e.g. 2.4 \cite{Prasad}.

We let $\mathbf{Q}$ be the $K_v$-parabolic subgroup of
$\mathbf{H}$ with
$$\mathbf{Q}(K_v)=\{\, g\in \mathbf{H}(K_v) \mid g a_+^\infty= a_+^\infty \,\}$$
The proof of this lemma amounts to showing that
$$\mathbf{P^+}(K_v) \leq \mathbf{Q}(K_v)$$

Note that $g a_+^\infty =a_+^\infty$ if and only if
$$d(ga_+^ne \,,\, a_+^ne)=d(a_+^{-n}ga_+^ne\,,\, e)$$ is a bounded
sequence. Since distance from $e$ is a proper function, and
because the action of $\mathbf{H}(K_v)$ on $X$ is proper, we can
alternatively characterize $\mathbf{Q}(K_v)$ as the group
$$\{\, g\in \mathbf{H}(K_v) \mid (a_+^{-n} g a_+^n)_{n \in
\mathbb{N}}\se \mathbf{H}(K_v) \mbox{ is precompact} \,\}$$

Let $\mathfrak{u}$ and $\mathfrak{h}$ be the Lie algebras of
$\mathbf{R_u(P^+)}$ and $\mathbf{H}$ respectively. We denote the
set of positive roots given by our ordering of $\Phi _K$ in
Section 3.2 as $\Phi _K^+ \se \Phi _K$, and we write the set of
roots that are linear combinations of elements in $\Delta _K
-\alpha _0$ as $[\Delta _K -\alpha _0]$.

If $\mathbf{T}$ is the maximal $K$-split torus in $\mathbf{H}$
that was chosen to produce the roots $\Phi _K$, then our choice of
$\mathbf{P^+}$ from Section 3.2 implies that
$$\mathfrak{u} = \bigoplus _{\alpha \in \Phi _K^+ - [\Delta _K -\alpha _0]}
\{\, v \in \mathfrak{h} \mid \text{Ad}\,t (v)=\alpha (t) v \mbox{
for all } t\in \mathbf{T} \,\}$$ Note that $\Phi _K^+ - [\Delta _K
-\alpha _0]$ is exactly the subset of $\Phi _K$ consisting of sums
of the form $\sum _{\alpha _i \in \Delta _K} n_i \alpha _i$ with
$n _i \geq 0$ for all $i$ and $n_0 \geq 1$. By our definition of
$\mathbf{T_1}\leq \mathbf{T}$ as being contained in the kernel of
every root in $\Delta _K - \alpha _0 $, and as $a_+ \in
\mathbf{T_1}$, we can express $\mathfrak{u}$ as a finite direct
sum
$$\mathfrak{u}= \bigoplus _{n\geq 1} \{\, v \in \mathfrak{u} \mid
\text{Ad}\,a_+ (v)=n \alpha _0 (a_+) v \,\}$$

 Since
$|\alpha _ 0 (a _+)|_v >1$, we see that for any $u \in
\mathbf{R_u(P^+)}(K_v)$,
$$a_+^{-n} u a_+^n \to 1$$ as $n \to
\infty$. Hence, if $u \in \mathbf{R_u(P^+)}(K_v)$ and $z \in
\mathbf{Z_H(T_1)}(K_v)$, then
$$a_+^{-n} uz a_+^n=a_+^{-n} u a_+^n z\to z$$ In particular, the above sequence
is precompact. As a consequence, $uz$, and thus all of
$$\mathbf{P^+}(K_v) =\mathbf{R_u(P)}(K_v) \rtimes
\mathbf{Z_H(T_1)}(K_v)$$ is contained in $\mathbf{Q}(K_v)$.

\ep

With a little more effort, it can be shown that $a_+^\infty \in
\Pi _+ ^\infty - \partial \Pi _+ ^\infty$, but we will not need
this fact.

Now we know the direction of $a_+^\infty$. Our last point of
business in this section is to see that this direction leads away
from the orbit $\mathbf{H}(\si)
 e$. This argument is standard.

\bigskip \noindent \textbf{Lemma 3.4(b)} \emph{For any $r>0$,
there exists an $n\in \mathbb{N}$ such that $a_+^n  e \notin
X(r)$.}

\pru Choose any nontrivial $\gamma \in
\mathbf{R_u(P^+)}(\mathcal{O}_S)$. As in the proof of the
preceding lemma, $ a_+^{-n} \gamma a_+^n \to 1$.

From Theorem I.1.12 of \cite{Ra}, the sequence $(a_+^n)_{n\in
\mathbb{N}} \se \mathbf{H}(K_S)$ induces a sequence in the
quotient space $\mathbf{H}(\si) \backslash \mathbf{H}(K_S)$ that
is not contained in any compact set. The lemma follows.

\ep

\medskip \noindent \large \textbf{3.5 A blueprint at infinity} \normalsize \medskip

In this section we will construct a cycle inside $X^\infty$ in the
direction
      given by the sequence $(a_+^{-n}  e)_{n\in \mathbb{N}}$. This is the direction
      in $\Sigma$ that is opposite to $a^\infty _+$.
      In Section 3.7, translates of this cycle will be ``coned off"
       from points of the form $a_+^n  e$. Then,
  these cones will be intersected with $\mathbf{R_u(P^+)}(K_S)V$ to produce
 cycles in $X(r_0)$.

      We let $\Pi _-^\infty$ be the simplex opposite to
      $\Pi _+^\infty$ in the spherical apartment $\Sigma ^\infty$.
      This simplex is the unique maximal simplex in $X^\infty$
      that is fixed under the action of $\mathbf{P^-}(K_S)$, where
      $\mathbf{P^-}$ is the maximal proper $K$-parabolic subgroup
       of $\mathbf{H}$ that contains $\mathbf{Z_H(T_1)}$ and
      is opposite to $\mathbf{P^+}$.

      We let $\Delta _-^\infty$ be the simplicial star of $\Pi _-^\infty$ in
      the apartment $\Sigma ^\infty$. That is, $\Delta _- ^\infty$ is the
union of all simplices in $\Sigma ^\infty$ that contain $\Pi
_-^\infty$.

\smallskip \noindent \textbf{The description of a chain in the boundary.}
Let $\sigma$ be a codimension $1$ simplex in $\Sigma ^\infty$ that
is contained in the boundary of $\Delta _- ^\infty$. The geodesic
continuation of $\sigma$ in $\Sigma _\infty$ is a great sphere,
that is, the boundary of a closed simplicial hemisphere $R_\alpha
\se \Sigma ^\infty$ (called a \emph{root space}). Among the two
possible hemispheres, $R_\alpha$ and $R_{- \alpha}$, in $\Sigma
^\infty$ that contain $\sigma$ in their boundary (called
\emph{opposite root spaces}), we fix notation so that $R_{-
\alpha}$ contains $\Pi _- ^\infty$.

\bigskip \noindent \textbf{Lemma 3.5(a)} \emph{There exists a group
element $u_{- \alpha} \in \mathbf{R_u(P^-)}(K_S)$ fixing $R_{-
\alpha}$ pointwise and satisfying the condition}
$$ \Sigma ^\infty
\cap u_{- \alpha}  \Sigma ^\infty = R_{- \alpha}$$

\pru We may assume that $S$ consists of a single valuation $v$.
The general case follows from the definition of the spherical
join.

Let $\mathbf{Q}$ be the minimal $K_v$-parabolic subgroup of
$\mathbf{H}$ corresponding to the chamber containing $\sigma$ and
$\Pi _-^\infty$. Let $\Phi _{K_v}$ be the set of roots of
$\mathbf{H}$ with respect to the maximal $K_v$-split torus in
$\mathbf{A}$, let $\Phi _{K_v}^{\text{nd}} \se \Phi _{K_v}$ be the
set of nondivisible roots, and let $\Delta _{K_v} \se \Phi
_{K_v}^{\text{nd}}$ be the set of simple roots associated with our
choice of $\mathbf{Q}$.

As explained in 5.6 of \cite{Tits spherical}, there is a root $-
\alpha \in \Phi _{K_v}^{\text{nd}}$ such that any nontrivial
element $u _{-\alpha}$ of the root group $\mathbf{U}_{(-
\alpha)}(K_v) \leq \mathbf{H}(K_v)$ fixes $R_{- \alpha}$ pointwise
and satisfies $ \Sigma ^\infty \cap u_{- \alpha}  \Sigma ^\infty =
R_{ - \alpha}$. (A similar statment holds by replacing $- \alpha$
throughout with $\alpha$, where the root $\alpha \in \Phi
_{K_v}^{\text{nd}}$ is the negative of $\alpha$.)

Note that all we have left to show is $\mathbf{U}_{(-\alpha)} \leq
\mathbf{R_u(P^-)}$.

Recall the standard correspondence that assigns to any subset $I
\se \Delta _{K_v}$ a $K_v$-parabolic subgroup of $\mathbf{H}$
containing $\mathbf{Q}$, denoted $\mathbf{Q}_I$; see e.g. 21.12
\cite{Borel text}. Since $\sigma $ is of codimension $1$ in
$\Sigma ^\infty$, the $K_v$-parabolic subgroup of $\mathbf{H}$
corresponding to $\sigma$ is of the form $\mathbf{Q}_{\{\beta\}}$
for a single simple root $\beta \in \Delta _{K_v}$. We also have
that $\mathbf{U}_{(-\alpha)}(K_v) \leq \mathbf{Q}_{\{\beta
\}}(K_v)$ and $\mathbf{U}_{(\alpha)}(K_v) \leq \mathbf{Q}_{\{\beta
\}}(K_v)$ since $\sigma \se R_{-\alpha} \cap R_{\alpha}$ is fixed
by $\mathbf{U}_{(-\alpha)}(K_v)$ and $\mathbf{U}_{(\alpha)}(K_v)$.
It follows from 21.12 of \cite{Borel text} that either $-\alpha
=\beta$ or $\alpha =\beta$.

Since $\mathbf{U}_{(-\alpha)}(K_v)$ fixes $R_{-\alpha}$ pointwise,
the chamber corresponding to $\mathbf{Q}$ is also fixed under the
action of $\mathbf{U}_{(-\alpha)}(K_v)$. Hence,
$\mathbf{U}_{(-\alpha)} \leq \mathbf{Q}$ implying that $-\alpha$
is positive under the ordering on $\Phi _{K_v}$ consistent with
$\Delta _{K_v}$. Now it must be that $-\alpha =\beta$.

Since $\Pi _- ^\infty \nsubseteq \sigma$, we have
$\mathbf{Q}_{\{-\alpha\}}=\mathbf{Q}_{\{\beta\}} \nleq
\mathbf{P^-}$. Therefore, if we assume $J \se \Delta _{K_v}$ is
such that $\mathbf{Q}_J=\mathbf{P^-}$, then $-\alpha \notin J$. It
follows that $\mathbf{U}_{(-\alpha)} \leq \mathbf{R_u(P^-)}$ as
desired.

\ep

Any $K_v$-parabolic subgroup of $\mathbf{H}$ that is contained in
$\mathbf{P^-}$ must contain $\mathbf{R_u(P^-)}$. Thus,
$u_{-\alpha} \in \mathbf{R_u(P^-)}(K_S)$ fixes $\Delta _-^\infty$
pointwise. Therefore, $\Delta _- ^\infty \se R_{- \alpha}$ which,
in turn, implies that $\Delta _- ^\infty \cap R_{\alpha }$ is the
union of some codimension $1$ simplices in the boundary of $\Delta
_- ^\infty$ (including $\sigma$). We name this union $F_\alpha$
and call it a \emph{geodesically continued face of} $\Delta _-
^\infty$. We take a minimal (hence finite) family of root spaces
$\{R_{\alpha }\}_{\alpha \in A}$ which exhaust the boundary of
$\Delta _- ^\infty$ as the union of the corresponding geodesically
continued faces of $\Delta _- ^\infty$.

Applying Lemma 3.5(a) to the opposite parabolic and opposite root
space, we have that for each for each $\alpha \in A$, there is a
group element $u_{\alpha } \in \mathbf{R_u(P^+)}(K_S)$ that fixes
$R_{\alpha }$ pointwise and satisfies $ \Sigma ^\infty \cap
u_{\alpha }  \Sigma ^\infty =R_{\alpha }$. Hence, $ \Delta _-
^\infty \cap u_{\alpha }  \Delta _- ^\infty = F_\alpha$.

We define the group $U \leq \mathbf{R_u(P^+)}(K_S)$ to be
generated by the finite set of $u_{\alpha }$ as above. As it will
be useful in Section 3.7, we also choose our $u_\alpha$ to fix the
point $e$. This can always be arranged by replacing the $u_\alpha$
with conjugates by elements of $\mathbf{A}(K_S)$.

It is well known that every element of $\mathbf{R_u(P^+)}(K_S)$
has order a power of $p$ (see e.g. 4.1 \cite{Borel text}), so $U$
is a $p$-group. Generalizing Schur's work on the generalized
Burnside problem, Kaplansky showed that any finitely generated
linear torsion group is finite (see e.g. Theorem 9.9 \cite{La}).
We conclude that $U$ is a finite $p$-group.

By abuse of notation, we shall denote the formal sum of chambers
in $\Delta _-^\infty$ simply by $\Delta _- ^\infty$. Now we form
the ($k(\mathbf{H},S)-1$)-chain $\sum _{u \in U} u  \Delta _-
^\infty $.

\smallskip \noindent \textbf{Properties of the chain in the boundary.}
 In the remainder of this section,
 we will show that $\sum _{u \in U} u  \Delta _-
^\infty $ is a cycle describing a simplicial decomposition of $ U
 \Delta _- ^\infty = \bigcup_{u \in U} u  \Delta _-
^\infty $.

\bigskip \noindent \textbf{Lemma 3.5(b)} \emph{If $u \in U$ is
nontrivial and $\mathfrak{C}^\infty \se \Delta _- ^\infty$ is a
chamber, then $u  \mathfrak{C}^\infty \nsubseteq \Delta _-
^\infty$.}

\pru Suppose $u  \mathfrak{C}^\infty \subseteq \Delta _- ^\infty$.
Then we have $\Pi _- ^\infty \se \mathfrak{C}^\infty \cap u
 \mathfrak{C}^\infty$ by the definition of $\Delta _-
^\infty$. As the action of $\mathbf{H}(K_S)$ on $ X_\infty$ is
type preserving, $u  \Pi _-^\infty = \Pi _-^\infty$. This implies
that $u \in
      \mathbf{P^-}(K_S) \cap \mathbf{R_u(P^+)}(K_S)=1$.

\ep

    \bigskip \noindent \textbf{Lemma 3.5(c)} \emph{The chain $\sum _{u \in
U} u  \Delta _- ^\infty $ is a
      cycle over $\mathbb{F}_p$.}

    \pru Suppose that $u \in U$ is nontrivial and that
      $  \Delta _- ^\infty \cap u  \Delta _- ^\infty$ contains an interior point $x$
      of a maximal simplex of a geodesically continued face of $\Delta _- ^\infty$, say $F_\alpha$.
      We begin by verifying that $u$ fixes $F_\alpha$ pointwise, and
     that $ F_\alpha =  \Delta _- ^\infty \cap u  \Delta _- ^\infty$.

      Indeed, $u$ fixes pointwise a simplex of
       $F_\alpha$ that contains $x$, since $u$ acts by type preserving simplicial
       automorphisms on $X ^ \infty$. The antipodal point of $x$ in $\Sigma ^\infty$
        is contained in the boundary of a chamber of $\Sigma ^\infty$ containing $\Pi
        ^\infty _+$; we call this chamber $\mathfrak{C}^\infty$. As in the comment immediately
        following proof of Lemma
        3.5(a), we see that $\mathfrak{C}^\infty$
         is fixed by $u \in \mathbf{R_u(P^+)}(K_S)$.

      The hemisphere $R _{\alpha }$ is the convex
      hull spanned by the simplex of $F_\alpha$ that contains $x$
       and the chamber $\mathfrak{C}^\infty$.
       Therefore, $u$
       fixes
      every point in $R_{\alpha } \supseteq F_\alpha$.

      For the remaining claim that
      $ F_\alpha =  \Delta _- ^\infty \cap u  \Delta _- ^\infty$:
       If there was a point $y\in  \Delta _- ^\infty \cap u
 \Delta _- ^\infty$
     outside of $R_{ \alpha}$, then
      $u$ would have to fix $y$ since $R_\alpha$ is fixed
      pointwise by $u$ and the action is by isometries. Hence, $u$ fixes pointwise
      the convex hull of $R_{ \alpha}$ and $y$. But
      that is all of
      $\Sigma ^\infty$, and any $u \in
      \mathbf{R_u(P^+)}(K_S)$ fixing $\Sigma ^\infty$ pointwise
      is the identity.
      So we have verified our claims.

We are now prepared to show that the homological boundary of $\sum
_{u \in U} u  \Delta _- ^\infty $  is $0$ modulo $p$. Applying the
boundary homomorphism yields:
$$\partial \Big(\sum _{u\in U} u  \Delta _- ^\infty \Big)
=\sum _{u \in U} \partial (u  \Delta _- ^\infty) =\sum_{u \in
U}\sum_{\alpha \in A} u  F_\alpha$$ where, again stretching
notation slightly, $F_\alpha$ denotes the formal sum of all
simplices in the geodesically continued face.

The claims we verified above show that, for $u,v \in U$ and all
$\alpha \in A$, either $u  F_\alpha \cap v  F_\alpha$ is contained
in the topological boundary of $u  F_\alpha$ or alternatively, $u
F_\alpha$ and $v  F_\alpha$ are equal as chains. Thus, we choose a
complete set $\{f_1, f_2, \ldots , f_n \}$ of representatives for
the chains in $\{u  F_\alpha\}_{u\in U,\alpha\in A}$ so that
$$\sum_{u \in U}\sum_{\alpha \in A} u  F_\alpha = \sum_{i=1}^n|U_{i}|f_i$$
where $U_i \leq U$ is the stabilizer of $f_i$. Since $U$ is a
finite $p$-group, $|U_i|$ is a power of $p$. Moreover, since each
$F_\alpha$ is stabilized by a nontrivial $u_\alpha \in U$, each
group $U_i$ is nontrivial. Therefore,
$$\partial \Big(\sum _{u\in U} u  \Delta _- ^\infty
\Big)\equiv  0 \text{ (mod $p$)}$$

    \ep

\noindent \textbf{Observation.} By the preceding lemmas, $U
 \Delta _- ^\infty$ represents a class in the homology group $\HOM
_{k(\mathbf{H},S)-1}( U  \Delta _- ^\infty \,;\,\mathbb{F}_p)$.

\medskip \noindent \large \textbf{3.6 A line of
 communication from infinity to $\mathbf{X(r_0)}$} \normalsize \medskip

In the next section, we will build cycles in $X(r_0)$ by
transferring the topological data from $U\Delta_-^\infty$ into
$X(r_0)$ by method of ``casting shadows" of $U\Delta _-^\infty$ on
$\mathbf{R_u(P^+)}(K_S)V$. For the shadow to contain the same
topological data as $U\Delta _-^\infty$, it is important, for
example, to have the shadow of $\Delta _-^\infty$ in $V$ be
compact. The purpose of this section is to establish that fact,
although we state this problem below using different language.

Recall that $\Sigma ^\infty$ can be regarded as the space of all
geodesic rays in $\Sigma$ based at $e$. We let $V^\infty \se
\Sigma ^\infty$ be the set of all geodesic rays contained in $V$
emanating from $e$. Note that $V^\infty$ is an equatorial sphere
in $\Sigma ^\infty$.

      We call a point in $\Sigma ^\infty$ \emph{rational} if it is
      represented by a geodesic ray based at $e$ that passes through
      another (and hence infinitely many) points of
$\mathbf{A}(K_S) e$. Let $\Sigma _\mathbb{Q}^\infty$
       denote the set of rational points in
      $\Sigma ^\infty$. Since $\mathbf{A}(K_S)$
      acts on $\Sigma$ as a lattice of translations of full
      rank $k(\mathbf{H},S)=\text{dim}(\Sigma)$, the set $\Sigma _\mathbb{Q}^\infty$
       is dense in $\Sigma ^\infty$.

      Similarly, we let $V_\mathbb{Q}^\infty$ denote the set of those
      points in $V^\infty$ that can be joined to $e$ by
      a geodesic ray passing through infinitely many points
      of
$\mathbf{A}(\si) e$. From our choice of $V$ before the Example in
Section 3.2, it is also clear that $V_{\mathbb{Q}}^\infty$ is
dense in $V^\infty$.

  \bigskip \noindent \textbf{Lemma 3.6(a)} \emph{We have $V_\mathbb{Q}^\infty = V^\infty
   \cap \Sigma _\mathbb{Q} ^\infty$.}

\pru The action of $\mathbf{A}(\si)$
        factors through the inclusion
$\mathbf{A}(\si) \hookrightarrow \mathbf{A}(K_S)$.
        Since $\mathbf{A}(\si)$
        acts on $V$ as a lattice of maximum rank
        $k(\mathbf{H},S)-1 = \text{dim}(V)$,
        the affine lattices $\mathbf{A}(\si)  e$
        and $V \cap \mathbf{A}(K_S)  e$
        are commensurable. Hence, they define identical rational
        structures at infinity.

      \ep

The goal of this section is:

\bigskip \noindent \textbf{Lemma 3.6(b)} \emph{We have $\Delta _- ^\infty
 \cap V ^\infty = \emptyset$.}

\pru We proceed by contradiction. So assume
         $ \Delta _- ^\infty \cap V^\infty \neq \emptyset$.
        Our first step will be to show that
        $ \Delta _- ^\infty \cap V _\mathbb{Q} ^\infty \neq \emptyset$.
        There are two cases. First, $V^\infty$ contains an interior point
        of $\Delta _- ^\infty$. Then the intersection
        $\Delta _- ^\infty \cap V^\infty$ is open in $V^\infty$ and
        contains a rational point since these are dense in
        $V ^ \infty$. That is $ \Delta _- ^\infty \cap V _\mathbb{Q} ^\infty \neq \emptyset$.
         Second, $V ^\infty$ contains a boundary simplex of $\Delta _- ^\infty$.
          Since the affine lattice
        $\mathbf{A}(K_S) e \se \Sigma$
        is commensurable to the affine lattice of vertices in the Euclidean
        Coxeter complex underlying the apartment $\Sigma$,
        rational points are dense in every simplex in $\Sigma ^\infty$.
        Therefore, $V^\infty \cap \big(\Sigma _\mathbb{Q} ^\infty
         \cap \Delta _- ^\infty \big) \neq \emptyset$.
        Using Lemma 3.6(a), we again find a point in
        $\Delta _- ^\infty \cap V _\mathbb{Q} ^\infty$.

        Now choose $b \in \mathbf{A}(\si)$ such that
        $b^n  e$ converges to a point
$b^\infty \in \Delta _- ^\infty \cap V _\mathbb{Q} ^\infty$
        as $n \to \infty$.

        Recall that for each $v\in S$, the group $\mathbf{R_u(P^-)}$ is
        contained in any minimal $K_v$-parabolic subgroup of $\mathbf{H}$ that is
        contained in $\mathbf{P^-}$. Therefore, $
        \mathbf{R_u(P^-)}(K_S)$ fixes $\Delta _-^\infty$ pointwise and, consequently, fixes
         the point $b^\infty \in \Delta _-^\infty$. As in
        the proof of Lemma 3.4(a),
        \begin{align*}
        \mathbf{R_u(P^-)}(\si ) & \leq
        \mathbf{R_u(P^-)}(K_S) \\
& \leq \{\, g\in \mathbf{H}(K_S) \mid g b^\infty=b^\infty \,\} \\
    & = \{\, g\in \mathbf{H}(K_S) \mid (b^{-n} g b^n)_{n \in
\mathbb{N}} \mbox{ is precompact} \,\}
\end{align*}
Therefore, for any $\gamma \in \mathbf{R_u(P^-)}(\si)$, the
sequence $(b^{-n} \gamma b^n)_{n \in \mathbb{N}} \se
\mathbf{H}(\si)$ is both discrete and precompact. Hence, it is
finite. We conclude that $$b^{-n}\gamma b^n = b^{-m}\gamma b^m$$
        for distinct $n$ and $m$. Now,
        $\gamma$ centralizes $b^{n-m}$.

        Let $\mathbf{D}$ be the subgroup of $\mathbf{A}$ that is the
        Zariski closure of the group generated by $b^{n-m}$. Then, $\gamma$ centralizes
        $\mathbf{D}$. Note that $b \in \mathbf{D}$, so $\mathbf{D}$ is also the
        Zariski closure of the group generated by $b$. Thus,
        $\mathbf{D}$ is independent of our choice of $\gamma \in
        \mathbf{R_u(P^-)}(\si)$. Hence, $ \mathbf{R_u(P^-)}(\si)$
        centralizes $\mathbf{D}$.

        Since $b^n e \to b^\infty$, iterates of $b$ define an unbounded sequence in
        $\mathbf{D}(K_v)$ for at least one $v \in S$. It follows that $\mathbf{D}$ contains a
        nontrivial $K_v$-split torus $\mathbf{D_d}$. Indeed, if $\mathbf{D}$ were
         $K_v$-anisotropic, then
         $\mathbf{D}(K_v)$ would be compact.

We denote the centralizer of $\mathbf{D_d}$ in $\mathbf{H}$ by
$\mathbf{L}$. Therefore, $\mathbf{L}$ is a Levi subgroup of a
$K_v$-parabolic subgroup of $\mathbf{H}$ (20.4 \cite{Borel text}).
It is clear that $\mathbf{A} \leq \mathbf{L}$.

We have shown that $$\mathbf{R_u(P^-)}(\si) \leq \mathbf{L}(K_v)$$
As $\mathbf{R_u(P^-)}$ is $K$-isomorphic as a variety to affine
space (see 21.20 \cite{Borel text}), $\mathbf{R_u(P^-)}(\si)$ is
Zariski dense in $\mathbf{R_u(P^-)}$ (use 3.1.1.ii \cite{Margulis
text}). Thus, $$\mathbf{R_u(P^-)} \leq \mathbf{L}$$

Since $\mathbf{L}(K_v)$ fixes $b^\infty \in \Delta ^\infty _-$,
there is a minimal $K_v$-parabolic subgroup of $\mathbf{H}$, which
we will write as $\mathbf{Q}$, such that $\mathbf{A} \leq
\mathbf{Q} \leq \mathbf{P^-}$ and such that $\mathbf{L}$ is a Levi
subgroup for a $K_v$-parabolic subgroup containing $\mathbf{Q}$.

Let $\Phi _{K_v}$ be the set of roots of $\mathbf{H}$ with respect
to a maximal $K_v$-split torus in $\mathbf{A}$, and let $\Phi
^+_{K_v}$ and $\Delta _{K_v}$ be the sets of positive and simple
roots respectively that correspond to $\mathbf{Q}$.

In the notation of 21.11 \cite{Borel text}, there is a proper
subset of simple roots $I \se \Delta_{K_v}$ such that
$$\mathbf{U}_{\Phi ^+_{K_v} - [I]} = \mathbf{R_u(P^-)}$$ where $[I]$
 is the subset of $\Phi _{K_v}$ consisting of
linear combinations of elements in $I$. There is also a subset of
simple roots $J \se \Delta _K$ such that
$$[J]=\{\, \alpha \in \Phi _{K_v} \mid \mathbf{U}_{(\alpha)} \leq
\mathbf{L} \,\}$$ If $\mathbf{L}$ is a proper subgroup of
$\mathbf{H}$, then $J$  is a proper subset of $\Delta _{K_v}$.

Since $\mathbf{H}$ is absolutely almost simple, $\Phi _{K_v}$ is
irreducible. Hence, there exists a ``highest root"
$\widetilde{\alpha} \in \Phi ^+ _{K_v}$ such that
$\widetilde{\alpha} = \sum _{\alpha \in \Delta _{K_v}} n_\alpha
\alpha$ with $n_\alpha \geq 1$ for all $\alpha \in \Delta _{K_v}$;
see e.g. VI.1.8 \cite{Bourbaki}. Thus, $\widetilde{\alpha} \in
\Phi ^+_{K_v} - [I]$ so
$$\mathbf{U}_{(\widetilde{\alpha})} \leq \mathbf{U}_{\Phi ^+_{K_v} - [I]}$$
 Stringing together the group equalities and
inclusions we have collected so far yields
$$\mathbf{U}_{(\widetilde{\alpha})} \leq \mathbf{U}_{\Phi ^+_{K_v}
- [I]} = \mathbf{R_u(P)} \leq \mathbf{L}$$

It follows from our description of the relationship between $[J]$
and the root groups contained in $\mathbf{L}$ that
$\widetilde{\alpha} \in [J]$. Hence, $J$ is not a proper subset of
$\Delta_{K_v}$, and we are left with $\mathbf{L}=\mathbf{H}$. That
is, the center of $\mathbf{H}$ contains the infinite group
generated by $b$. This is our contradiction.

 \ep

\medskip \noindent \large \textbf{3.7 Cycle assembly in $\mathbf{X(r_0)}$} \normalsize \medskip

Let $C_\Delta \se \Sigma $ be the cone of all geodesic rays
contained in $\Sigma$, based at $e$, and limiting to points in
$\Delta _- ^\infty$.

Recall our choice of $a_+ \in \mathbf{H}(K_S)$ as a translation of
$\Sigma $ such that $$a_+^n e \to a_+^\infty \in \Pi _+^\infty $$
 Recall also that $\Pi ^\infty _+ \se \Sigma ^\infty$ is the
collection of antipodal points for points in $\Pi ^\infty _- \se
\Sigma ^\infty$ and that $\Delta _-^\infty$ is the union of
chambers in $\Sigma ^\infty $ containing $\Pi ^\infty _-$.

Therefore, Lemma 3.6(b) implies for $n \geq 1$ that any geodesic
ray emanating from $a_+^n e$ and limiting to $\Delta _- ^\infty$
is separated by $V$. Hence, there is a well-defined geodesic
projection toward $a_+^n e$ that gives rise to a homeomorphism
$$\Delta _-^\infty \longrightarrow V \cap a_+^n  C_\Delta$$

Recall that we chose $U$ to fix $e$. Thus, $a_+^n U a_+^{-n}$
fixes $a_+^n  e$. It follows that for all $u \in U$, there are
well-defined geodesic projections toward $a_+^n  e$ that give rise
to homeomorphisms
$$a_+^n u  \Delta _- ^\infty = a_+^n u a_+^{-n}  \Delta _- ^\infty
\longrightarrow a_+^n u a_+^{-n}  (V \cap a_+^n  C_\Delta )$$ Note
that these maps piece together to give a continuous surjection
$$\pi _n : a_+^n U  \Delta _- ^\infty \longrightarrow a_+^n U a_+^{-n}
  (V \cap a_+^n  C_\Delta )$$
whose image is contained in $\mathbf{R_u(P^+)}(K_S)  V \se
X(r_0)$. The collection of $\sum _{u \in  U} \pi _n \big( a_+^n  u
\Delta _- ^\infty  \big)$ are the cycles we have been searching
for throughout this paper.

\bigskip \noindent \textbf{Lemma 3.7(a)} \emph{There is a point $s
\in \Sigma$, a chamber $\mathfrak{s} \se \Sigma$, and a sector
$\mathfrak{S} \se C _\Delta$ such that:}
\begin{itemize}
        \item[(\emph{i})]
          $s \in \mathfrak{s} \se \mathfrak{S}$; \emph{and}
        \item[(\emph{ii})]
           \emph{For each nontrivial $u \in U$,}  $$\mathfrak{S} \cap \varrho _{\Sigma ,  \mathfrak{s}}
          ( u  C_\Delta) = \emptyset$$ \emph{where $\varrho _{\Sigma ,  \mathfrak{s}} : X \rt
          \Sigma$ is the building retraction for the pair
          $(\Sigma , \mathfrak{s})$.}
      \end{itemize}

\pru Let $\mathfrak{S}' \se C_\Delta$ and $\mathfrak{T}\se X $ be
sectors that do not contain a common subsector. Consider an
apartment $\Sigma _* \se X$ that contains disjoint subsectors
        $\mathfrak{S} _0 \se \mathfrak{S}'$ and
        $\mathfrak{T} _0 \se \mathfrak{T}$. For any chamber
        $\mathfrak{c} \se \mathfrak{S}_0$, the retraction
         $\varrho _{\Sigma , \mathfrak{c}}$ restricts to an
         isometry from $\Sigma _*$ to $\Sigma$ that fixes
         $\mathfrak{S}_0$ pointwise. Thus, we have $$\mathfrak{S}_0
\cap \varrho _{\Sigma , \mathfrak{c}}(\mathfrak{T}_0)=\emptyset$$

Choose $D \geq 0$ such that $\mathfrak{T}$ is contained within the
closed metric $D$-neighborhood of $\mathfrak{T}_0$. Now choose
$\mathfrak{S} \se \mathfrak{S}_0$ such that the closed metric
$D$-neighborhood of $\mathfrak{S}$ in $\Sigma$ is completely
contained within $\mathfrak{S}_0$. Then for any two chambers
$\mathfrak{s} \se \mathfrak{S}$ and $\mathfrak{t} \se
\mathfrak{T}$, the distance from $\varrho _{\Sigma ,
\mathfrak{s}}(\mathfrak{t})$ to $\varrho _{\Sigma , \mathfrak{s}}(
\mathfrak{T}_0)$ is at most $D$ since $\varrho _{\Sigma ,
\mathfrak{s}}$ does not increase distances. As the distance from
$\varrho _{\Sigma , \mathfrak{s}}(\mathfrak{T_0})$ to
$\mathfrak{S}$ is at least $D$, we find $$\mathfrak{S} \cap
\varrho _{\Sigma , \mathfrak{s}}(\mathfrak{T})=\emptyset $$

By Lemma 3.5(b), $u  C_\Delta $ can be covered by finitely many
$\mathfrak{T}$ as above for any nontrivial $u \in U$. Thus, we can
assume, after perhaps passing to a subsector of $\mathfrak{S}$,
that $\mathfrak{S} \cap \varrho _{\Sigma ,
\mathfrak{s}}(\mathfrak{T})=\emptyset $ for all such $u$ and
$\mathfrak{T}$. Hence, the lemma is satisfied for any choice of $s
\in \mathfrak{s}$.

\ep

We fix $s$, $\mathfrak{s}$, and $\mathfrak{S}$ as above, and for
every $n \in \mathbb{N}$ we let $$\phi _n : \Sigma - \{a_+^n
 s \} \rt \Sigma ^\infty$$ be the visual projection to the
boundary from the point $a^n _+ s$.

\bigskip \noindent \textbf{Lemma 3.7(b)} \emph{For every $r \geq r_0$,
the inclusion $X(r_0) \hookrightarrow X(r)$ induces a nontrivial
homomorphism
$$\HOM_{k(\mathbf{H},S)-1}(X(r_0)\,;\,\mathbb{F}_p) \longrightarrow
\HOM_{k(\mathbf{H},S)-1}(X(r)\,;\,\mathbb{F}_p)$$}

\pru Choose $n \in \mathbb{N}$ such that $V$ separates $a_+^n
 \mathfrak{S} $ into a compact component (containing $s$)
and a noncompact component, and such that $$a_+^n  e \notin
X\big(r +d(e,s)\big)$$ The latter condition can be arranged by
Lemma 3.4(b), and it implies that $$\varrho _{\Sigma ,\, a_+^n
\mathfrak{s}} ^{-1}(a^n_+ s)=\{a_+^n   s \} \nsubseteq X(r)$$
where $\varrho _{\Sigma ,\, a_+^n  \mathfrak{s}}$ is the
retraction corresponding to the pair $(\Sigma ,\, a_+^n
\mathfrak{s})$ Therefore, the following composition is well
defined:
$$  a_+^n U  \Delta _-
^\infty  \rightarrow X(r_0) \hookrightarrow X(r) \rt \Sigma \rt
\Sigma ^\infty$$ where the map on the left is $\pi _n$, the map
second from the right is $\varrho _{\Sigma ,\, a_+^n
\mathfrak{s}}$, and the map on the far right is $\phi _n$.

Since $\varrho _{\Sigma ,\, a_+^n  \mathfrak{s}}$ is simply
$\varrho _{\Sigma , \mathfrak{s}}$ conjugated by $a_+^n$, Lemma
3.7(a) implies that there is an open neighborhood of $\phi _n (V
\cap a_+^n  \mathfrak{S} ) \se \Delta_- ^\infty$ that has
$1$-point pre-images of points under the above composition. Hence,
using excision---as in determining degrees of maps between spheres
(see e.g. Proposition 2.30 of \cite{Hatcher})---one sees that the
induced homomorphism
$$\HOM_{k(\mathbf{H},S)-1}(a_+^n U \Delta _-
^\infty\,;\,\mathbb{F}_p) \longrightarrow
\HOM_{k(\mathbf{H},S)-1}(\Sigma ^\infty \,;\,\mathbb{F}_p)$$ is
nontrivial. Our result follows as the above homomorphism factors
through
$$\HOM_{k(\mathbf{H},S)-1}(X(r_0)\,;\,\mathbb{F}_p) \longrightarrow
\HOM_{k(\mathbf{H},S)-1}(X(r)\,;\,\mathbb{F}_p)$$

\ep

Our proof of Proposition B* is complete.

\medskip \noindent \large \textbf{4 Appendix} \normalsize
\medskip

For completeness, in this section we record two results.

\smallskip \noindent \textbf{Existence of anisotropic tori.}
 The following lemma is well-known, though we do not know of a reference for it.
 We thank Stephen DeBacker for communicating the following proof to
 us.

\bigskip \noindent \textbf{Lemma 4(a)} \emph{Any semisimple group defined over
        $\mathbb{F}_q((t))$ contains
        a maximal
        $\mathbb{F}_q((t))$-torus that is
        $\mathbb{F}_q((t))$-anisotropic.}

\pru Let $k= \mathbb{F}_q((t))$, let $k_s$ be a separable closure
of $k$, and let $\overline{\mathbb{F}}_q$ be the algebraic closure
of $\mathbb{F}_q$ in $k_s$.

We define $\tau_q :\overline{\mathbb{F}}_q \rt
\overline{\mathbb{F}}_q$ as the Frobenius automorphism $x \mapsto
x^q$, and we extend $\tau _q$ to an automorphism
$\widehat{\tau}_q: \overline{\mathbb{F}}_q((t)) \rt
\overline{\mathbb{F}}_q((t))$ by applying $\tau_q$ to the
coefficients of the Laurent series. Finally, we choose some $\tau$
in the Galois group $\text{Gal}(k_s/k)$ that extends
$\widehat{\tau}_q$. Below we will write $k_u$ for the field
$\overline{\mathbb{F}}_q((t))$ so that we have $k<k_u<k_s$.

        Suppose $\mathbf{G}$ is a connected semisimple
       $k$-group, and let $\mathbf{T} \leq \mathbf{G}$
             be a maximal $k_u$-split torus that is defined over
             $k$. Such a torus exists by
        a result of Bruhat-Tits (Corollaire 5.1.12 of \cite{Bruhat-Tits}).
        We denote the character group of $\mathbf{T}$ by $\mathbf{X}(\mathbf{T})$.

 Let $\mathbf{N}$ be the normalizer of $\mathbf{T}$ in
$\mathbf{G}$, and let $\mathbf{Z}$ be the centralizer of
$\mathbf{T}$ in $\mathbf{G}$. Note that $\mathbf{Z}$ is defined
over $k$, and $\mathbf{Z}$ is a torus since $\mathbf{G}$ is
necessarily $k_u$-quasisplit (see e.g. Lemma 4.7 \cite{Sp 2}).
Choose an $n \in \mathbf{N}(k_u)$ such that $n \circ \tau$ acts as
a so-called ``twisted Coxeter element"  on the vector space
$\mathbf{X}(\mathbf{T}) \otimes \mathbb{R}$ (see e.g. Lemma 7.4(i)
\cite{Sp}); thus, $n \circ \tau$ has no nontrivial fixed points in
$\mathbf{X}(\mathbf{T}) \otimes \mathbb{R}$.

By Lang's theorem, there exists a $g \in \mathbf{G}(k_u)$ such
that $g^{-1}\tau (g) =n$ (see e.g. the proof of 16.4 \cite{Borel
text}). We claim that $\sideset{^g}{}{\cent}$ satisfies the lemma.

First, to show that $\sideset{^g}{}{\cent}$ is defined over $k$,
we need to verify that $\sideset{(^g}{)^\varphi}{\cent}  =
\sideset{^g}{}{\cent}
 $ for all $\varphi \in
\text{Gal}(k_s/k)$. Since $\sideset{^g}{}{\cent}$ is defined over
$k_u$, it suffices to check that $\sideset{(^g}{)^\tau}{\cent} =
\sideset{^g}{}{\cent}
 $. Thus we see that
$\sideset{^g}{}{\cent}$ is defined over $k$:
$$(\sideset{^g}{}{\cent})^\tau = \sideset{^{\tau (g)}(}{^\tau )}{\cent} =
\sideset{^{gn}}{}{\cent} =\sideset{^g}{}{\cent}$$

Second, note that showing $\sideset{^g}{}{\cent}$ is
$k$-anisotropic reduces to verifying that $\sideset{^g}{}{\torus}$
is $k$-anisotropic since the maximal $k_u$-anisotropic torus of
$\mathbf{Z}$ remains $k_u$-anisotropic (and hence $k$-anisotropic)
after conjugating by $g$.

To show that $\sideset{^g}{}{\torus}$ is $k$-anisotropic, it
suffices to check that there are no nontrivial fixed points of
$\tau $ on the character group $\mathbf{X}( \sideset{^g}{}{\torus}
)$. As shown in Proposition 3.3.4(i) of \cite{Carter}, there is a
bijection $\mathbf{X}( \sideset{^g}{}{\torus} ) \rt \mathbf{X}(
\mathbf{T} )$ that creates a correspondence between the $\tau
$-action on $\mathbf{X}( \sideset{^g}{}{\torus} )$ and the $(n
\circ \tau) $-action on $\mathbf{X}( \mathbf{T} )$. Hence, the
result follows from our choice of $n$.

\ep

\smallskip \noindent \textbf{The ``if" implication of Theorem A.}
We sketch the proof of this result as can be found in
Th\'{e}or\`{e}me 4 of \cite{Serre initial}.

We may assume
      that $\mathbf{G}$ is connected.
      The arithmetic group $\mathbf{G}(\si)$ is a direct product
      ``up to finite groups'' of $\gd (\si)$ and
      $\mathbf{Z(G)}(\si)$, where $\mathbf{Z(G)}$
      is the center of $\mathbf{G}$. By Dirichlet's units theorem,
      $\mathbf{Z(G)}(\si)$ is a finitely generated
      abelian group, hence of type $FP_\infty$.

      Since $\gd$ is $K$-anisotropic, Behr-Harder reduction theory
      yields that
      $\gd (\si)$ is a cocompact lattice in a semisimple group. Hence, $\gd (\si)$ acts properly and
      cocompactly on a Euclidean building. Therefore, $\gd (\si)$
      is of type $FP_\infty$ as well. We conclude that
      $\mathbf{G} (\si)$ is of type $FP_\infty$. (More is true: since $\gd (\si)$
       is residually finite, and there are, up to conjugacy, only finitely many
       stabilizer subgroups of $\gd (\si)$ for the action on the
       building,
      $\gd (\si)$ contains a finite-index torsion-free subgroup $\Gamma $
       that acts freely and cocompactly. Thus, there is a finite $K(\Gamma ,1)$.)

\bigskip

\noindent Kai-Uwe Bux
\newline \noindent Department of Mathematics
\newline \noindent P. O. Box 400137 (Kerchof Hall)
\newline \noindent University of Virginia
\newline \noindent Charlottesville, VA 22904-4137
\newline \noindent Email:  bux\_2004@kubux.net

\bigskip

\noindent Kevin Wortman
\newline \noindent Department of Mathematics
\newline \noindent Cornell University
\newline \noindent Malott Hall
\newline \noindent Ithaca, NY 14853
\newline \noindent Email: wortman@math.cornell.edu

\end{document}